\documentclass[smallextended,referee,envcountsect]{svjour3}
\smartqed
\usepackage{graphicx}
\usepackage[]{inputenc,amsmath,amssymb}

\DeclareMathOperator{\erf}{erf}
\usepackage[dvipsnames]{xcolor}
\usepackage{booktabs}   
\usepackage{todonotes}
\usepackage{mathtools}
\usepackage{multicol}
\journalname{JOTA}
\usepackage{lscape} 
\usepackage{soul}

\begin{document}

\title{Mixed Finite Differences Scheme for Gradient Approximation}

\author{Marco Boresta \and Tommaso Colombo \and Alberto De Santis \and Stefano Lucidi}

\institute{Department of Computer, Control and Management Engineering Antonio Ruberti, Sapienza
University of Rome, Via Ariosto 25, 00185 Rome, Italy \\
              \{boresta, colombo, desantis, lucidi\}@diag.uniroma1.it
}

\date{Received: date / Accepted: date}

\maketitle

\begin{abstract}

In this paper we focus on the linear functionals defining an approximate version of the gradient of a function.
These functionals are often used when dealing with optimization problems where the computation of the gradient of the objective function  is costly or the objective function values are affected by some noise.
These functionals have been recently considered to estimate the gradient of the objective function by the expected value of the function variations in the space of directions. The expected value is then approximated by a sample average over a proper (random) choice of sample directions in the domain of integration. In this way the approximation error is characterized by statistical properties of the sample average estimate, typically its variance. Therefore, while useful and attractive bounds for the error 
variance can be expressed in terms of the number of function evaluations, 
nothing can be said on the error of a single experiment that could be quite large. This work instead is aimed at deriving an approximation scheme for linear functionals approximating the gradient, whose  error of approximation can be characterized by a deterministic point of view in the case of noise-free data.
The previously mentioned linear functionals are no longer considered as expected values over the space of directions, but rather as the filtered derivative of the objective function by a Gaussian kernel.
By using this new approach, a gradient estimation based on a suitable linear combination of central finite differences at different step sizes is proposed and deterministic bounds that do not depend on the particular sample of points considered are computed. In the noisy setting, on the other end, the variance of the estimation error of the proposed method is showed to be strictly lower than the one of the estimation error of the Central Finite Difference scheme.
Numerical experiments on a set of test functions are encouraging, showing good performances compared to those of some methods commonly used in the literature, also in the noisy setting.
\end{abstract}


\section{Introduction}

DFO algorithms have become increasingly important since they provide a proper methodology to tackle most of the optimization problems considered in various fields of application. As reported in \cite{berahas2019theoretical,fazel2018global,salimans2017evolution} typical applications fall within the simulation-based optimization problems such as policy optimization 
in reinforcement learning.
DFO methods arise when derivative information is either unavailable, or quite costly to obtain, not to mention when only noisy sample of the objective function are available. In the latter case it is known that most methods based on finite difference is of little use  \cite{kolda2003optimization,wild2008orbit}.\hfill \break
\indent One of the approaches in DFO algorithms is that of computing a proper estimate of the gradient of the objective function. Finite difference approximation schemes were already present in early times \cite{polyak1987introduction} and have recently been reconsidered as sample average approximations of functionals defining a "filtered version" of the objective function \cite{flaxman2004online,nesterov2017random,balasubramanian2018zeroth,berahas2019linear}. These functionals arise when defining a gradient approximation as the average of the function variation along all the directions in the whole space.
In the most popular methods the average is performed by weighting the function variations along directions generated either with a uniform kernel on the unit ball \cite{flaxman2004online}, or with a Gaussian kernel \cite{balasubramanian2018zeroth}. These integrals are considered as ensemble averages over the space of the directions of differentiation, and then are approximated by sample averages over a random sample of directions, with various methods. As a general policy, the approximation error is then characterized by its statistical properties (even in the noise-free setting), the variance is expressed in terms of the number of function calculations, and nice bounds are provided to trade-off precision of the gradient estimation and computational costs. Nevertheless it is plain that the error on a single sample may be quite large, even though its variance is bounded.

In this paper we focus on a different point of view. The functional defining a filtered version of the objective function is considered as weak derivative of the objective function rather than expected values over the space of the directions \cite{ziemer2012weakly}. The gradient estimation is therefore obtained by considering a numerical approximation of the functional integral, and the estimation error is evaluated in a deterministic fashion. The estimate is obtained by a suitable linear combination of central finite differences at steps with increasing size.
Bounds on the approximation error with the proposed method are derived, and the variance of the error in the case of noisy data is also presented.

The goodness of the approximation is experimentally evaluated by comparing the proposed method with those considered benchmarks by the literature -  namely: Forward Finite Differences (FFD), Central Finite Differences (CFD) \cite{polyak1987introduction},  Gaussian Smoothed Gradient (GSG), Central Gaussian Smoothed Gradient (cGSG) \cite{nesterov2017random,flaxman2004online}
- over the benchmark of the Schittkowski functions \cite{schittkowski2012more}.  Encouraging results are obtained, both in the noise-free and in the noisy setting. 

\noindent
The paper is organized as follows: 
section \ref{The gradient estimation} formally introduces the gradient estimation problem, highlighting the difference between the approach proposed in this article and the one of several estimates proposed in the literature.
In section \ref{A new estimate of the gradient} we present the proposed approximation scheme - NMXFD, with an emphasis on its link with the Finite Difference Method. A theoretical comparison between the variance of the estimation errors of the proposed method and of the CFD scheme is proposed in section \ref{Noisy_data}.
Section \ref{Numerical experiment} presents numerical results and conclusions are drawn in section \ref{conclusioni}.

\section{The Gradient Estimate}
\label{The gradient estimation}
In this paper we consider the following unconstrained optimization problem in the derivative free optimization (DFO) setting \cite{larson_menickelly_wild_2019,conn2008geometry}:
\begin{equation}
    \min_{x \in R^n}  f(x)
    \label{problema_partenza}
\end{equation}
\noindent
where $f:\, R^n\mapsto R$ is a function with continuous derivative, i.e. $f\in\mathcal{C}^1(R^n)$, and we denote the gradient $\nabla f:\, R^n\mapsto R^n$ such that for any $x\in R^n$
$$
\nabla f(x)=\begin{bmatrix} \frac{\partial f}{\partial x_1}(x)\\ \vdots 
\\ \frac{\partial f}{\partial x_n}(x)\end{bmatrix}.
$$
In this section the problem of a numerical approximation of the gradient $\nabla f(x)$ is considered.
The most popular approximation scheme is the standard finite difference method  \cite{polyak1987introduction}, but interesting alternative schemes 
are proposed in papers \cite{flaxman2004online,balasubramanian2018zeroth}. A general estimate is obtained according to the following formula:
\begin{equation}
\label{formula2}
    G_{\sigma}(x) := \frac{1}{\sigma} \int_{R^n}^{} f(x + \sigma s) \,  
s \, \varphi(s) \,ds
\end{equation}
\noindent where $\varphi(s):\, R^n \mapsto R$ denotes either a standard Gaussian Kernel $\mathcal{N}(0, I_n )$ or a uniform kernel on the unit ball $\mathcal{B}(0, 1 )$, $ds=ds_1\cdot ds_2\cdot\ldots\cdot ds_n$ is the 
volume element in $R^n$, and $\sigma>0$ is a scale parameter.
The approximation error has different bounds depending on the assumptions on $f$ (see \cite{berahas2019theoretical}).
If the function $f$ is continuously differentiable, and its gradient is L-Lipschitz continuous for all $x \in R^n$, then
\begin{equation}
\label{bound_1}
    || G_{\sigma}(x) - \nabla f(x)  || \leq C_{\varphi} L \sigma
\end{equation}
\noindent
where $C_{\varphi}$ is a positive constant whose value depends on the kernel.
If the function $f$ is twice continuously differentiable, and its Hessian 
is H-Lipschitz continuous for all $x \in R^n$, then
\begin{equation}
\label{bound_2}
     || G_{\sigma}(x) - \nabla f(x)  || \leq C_{\varphi} H \sigma^2.
\end{equation}
\noindent Both bounds (\ref{bound_1}) and (\ref{bound_2}) show that
$$
\lim_{\sigma \to 0} G_{\sigma}(x) = \nabla f(x).
$$
We will now work out formula (\ref{formula2}) considering the (standard) Gaussian kernel
\begin{equation}
    \varphi(s) \sim \mathcal{N}(0, I_n ) = \frac{1}{(\sqrt{2 \pi})^n} \exp{\{{-\frac{1}{2}\sum_{i=1}^n {s_i}^2}\}} = \prod_{i=1}^n \varphi(s_i)
\end{equation}
\noindent
but the considerations that follow hold also if a uniform kernel over the unit 
ball is considered. \\
Let us consider this further notation: for any $x\in R^n$ denote by $\bar 
x_i\in R^{n-1}$ the following vector $ \begin{bmatrix} x_1, x_2 \, \ldots \, ,x_{i-1}, x_{i+1}, \, \ldots \, , x_n\end{bmatrix}^T$. With some abuse of notation, but for sake of simplicity in the use of formulas, when addressing a given coordinate $x_i$ in a vector $x$ let us write $x$ as $[x_i\>\bar x_i]^T$ and denote $f(x)$ as $f(x_i,\bar x_i)$ and $\varphi(s)=\varphi(s_i)\varphi(\bar s_i)$, with $\varphi(\bar s_i)=\prod_{j\ne i}^n \varphi(s_j)$; consistently, the volume element becomes $ds = ds_i\cdot d\bar s_i$. In case of a vector function $f(z)$, to address explicitly its $i-th$ entry we write it as $[(f(z))_i\>\overline{(f(z))}_i]^T$. Then, estimate (\ref{formula2}) is rewritten as follows
\begin{align}
    G_{\sigma}(x) = \frac{1}{\sigma}\int_{R^n} f(x_1 + \sigma s_1, \ldots, x_n + \sigma s_n) \begin{pmatrix} s_1\\\vdots\\s_n\end{pmatrix} \prod_{i=1}^n \varphi(s_i)\,   ds
\end{align}
\begin{gather}
 =
  \begin{bmatrix}
\frac{1}{\sigma}  \int_{R^n} f(x_1 + \sigma s_1, \bar x_1 + \sigma \bar s_1) \, s_1 \,\varphi(s_1)  \varphi(\bar s_1) \, ds_1 \, d\bar s_1\\
\vdots\\
 \frac{1}{\sigma}  \int_{R^n}f(x_i + \sigma s_i, \bar x_i + \sigma \bar s_i) \, s_i \, \varphi(s_i)  \varphi(\bar s_i) \, ds_i \, d\bar s_i\\
\vdots\\
 \frac{1}{\sigma}  \int_{R^n} f(x_n + \sigma s_n, \bar x_n + \sigma \bar s_n) \, s_n \, \varphi(s_n)  \varphi(\bar s_n) \, ds_n \, d\bar s_n
   \end{bmatrix}
   \label{vect}
\end{gather}
Let us consider the generic entry of vector (\ref{vect})
\begin{align}
    (G_{\sigma}(x))_i &= \frac{1}{\sigma} \int_{R^n} f(x_i + \sigma s_i, \bar x_i + \sigma \bar s_i) s_i \varphi(s_i)  \varphi(\bar s_i) ds_i\, d\bar s_i
 \label{entry}
\end{align}
By the Fubini theorem we can compute it as follows
\begin{equation}
(G_{\sigma}(x))_i =
\int_{R^{n-1}}  \varphi(\bar s_i) \Bigg( \frac{1}{\sigma} \int_{-\infty}^{+\infty} f(x_i + \sigma s_i, \bar x_i + \sigma \bar s_i) s_i \varphi(s_i)   ds_i    \Bigg )d\bar s_i
\label{Idpartialxi}
\end{equation}
The expression in parentheses is  the estimate of the directional derivative of $f(x)$ along the $i$-th coordinate $x_i$ and computed at the point 
 $(x_i, \bar x_i + \sigma \bar s_i)$, i.e.
\begin{equation}
g_\sigma(x_i, \bar x_i + \sigma \bar s_i)  :=
 \frac{1}{\sigma} \int_{-\infty}^{+\infty} f(x_i + \sigma s_i, \bar x_i + 
\sigma \bar s_i) s_i \varphi(s_i)   ds_i.
\label{partialxi}
\end{equation}
\noindent Hence expression (\ref{entry}) becomes
\begin{equation}
    (G_{\sigma}(x))_i = \frac{1}{\sigma} \int_{R^{n-1}} g_{\sigma}(x_i, \bar x_i 
+ \sigma \bar s_i) \, \varphi(\bar s_i) \,  d\bar s_i.
    \label{f6}
\end{equation}
Therefore, the generic entry of the gradient estimate $G_{\sigma}(x)$ in formula (\ref{vect}) is the average of function (\ref{partialxi}) weighted by a $(n-1)$-dimensional Gaussian kernel $\varphi(\bar s_i)=\mathcal{N}(0,  I_{n-1} )$ over the subspace $R^{n-1}$ of $R^n$. As a consequence, 
the computation of any entry of vector $G_{\sigma}(x)$ implies an integration over $R^n$.
In papers \cite{balasubramanian2018zeroth,berahas2019linear} this problem 
is overcome by considering that (\ref{formula2}) is indeed an ensemble average of function $f(x + \sigma s)  s$ over all the directions $s\in R^n$ weighted by the Gaussian distribution $\varphi(s)\sim\mathcal{N}(0, I_n )$. Therefore we can write
\begin{equation}
    G_{\sigma}(x) = \frac{1}{\sigma} E_{\varphi}[f(x + \sigma s)  s].
\end{equation}
Now the ensemble average can be well approximated by sampling a set of $M$ independent directions $\{s_i\}$  in $R^n$ according to $\mathcal{N}(0, I_n )$, and considering the \emph{sample average approximation} of 
$E_{\varphi}[f(x + \sigma s)  s]$
\begin{equation}
    G_{\sigma}(x) \simeq \frac{1}{M}\sum_{i=1}^M \frac{(f(x + \sigma s_i)-f(x))s_i}{\sigma}.
    \label{fd}
\end{equation}
or its simmetric version
\begin{equation}
    G_{\sigma}(x) \simeq \frac{1}{M}\sum_{i=1}^M \frac{(f(x + \sigma s_i)-f(x-\sigma s_i))s_i}{2\sigma}.
    \label{cfd}
\end{equation}
The same argument holds if a uniform distribution over the unit ball is considered for the ensemble average \cite{flaxman2004online}. Now, only $M+1$ function computations in case of (\ref{fd}) or $2M$ in case of (\ref{cfd}) are needed and the convergence properties of the sample estimate to the ensemble average are well established: the sample average is an unbiased estimate and its accuracy increases with increasing $M$. In \cite{berahas2019linear}, suitable expressions of the estimation error variance are 
found in terms of the number of samples $M$ and the values of  some smoothness parameters of function $f$. Therefore very useful formulas are given that define the required sample size to obtain a chosen accuracy, with a fixed level of confidence $1-\alpha$. This is a typical statistical characterization of the error, that is robust over the whole ensemble of possible trials, but of course leaves a risk $\alpha$ to have a large error on a single experiment.

In this paper, by exploiting formula (\ref{partialxi}), the following gradient estimate is proposed
\begin{equation}\overline G_\sigma(x) := \begin{bmatrix}
g_\sigma(x_1, \bar x_1), \hdots ,\, g_\sigma(x_i, \bar x_i), \hdots\, ,g_\sigma(x_n, \bar x_n)
\end{bmatrix}^T
\label{qgrad}
\end{equation}
where   \begin{equation}
\label{single_component}
    g_\sigma(x_i, \bar x_i) =  \frac{1}{\sigma} \int_{-\infty}^{+\infty} f(x_i + \sigma s_i, \bar x_i )\, s_i \,\varphi(s_i) \,  ds_i
\end{equation}
is obtained from (\ref{partialxi}) with $\bar s_i = 0 , \>  i = 1,..., d$.
This is a different result from estimate (\ref{vect}), and appears to be more practical since only line integrals are involved in the formula. \\
The following theorem shows that estimate $\overline G_\sigma(x)$ is close to $G_\sigma(x)$ and converges to it as $\sigma$ tends to zero.
\begin{theorem}
\label{prova_teorema0}
Let $\nabla f(x)$ be Lipschitz continuous with constant \textit{L} for all $x \in R^n$.
Then we have that
\begin{equation}
    || G_{\sigma}(x) - \overline G_\sigma(x) || \leq L\,\sigma\,\sqrt{n(15 + 7(n-1))}.
\end{equation}
\end{theorem}
\noindent
\textbf{Proof}: See Appendix for the proof.

Next theorem shows that $\overline G_\sigma(x)$ is indeed a good approximation of the true gradient $\nabla f(x)$ and converges to it as $\sigma$ tends to zero.
\begin{theorem}
\label{prova_teorema}
Let $f(x)$ be  continuously differentiable for all $x \in R^n$.
The following holds:
\begin{equation}
\label{tesi}
    \lim_{\sigma\to 0}  \overline G_\sigma(x) = \nabla f(x)
\end{equation}
\end{theorem}
\noindent
\textbf{Proof}:
We prove (\ref{tesi}) component-wise. By integration by parts we have
\begin{align}
    g_\sigma(x_i, \bar x_i) = &\frac{1}{\sigma} \int_{-\infty}^{+\infty} f(x_i + \sigma s_i, \bar x_i) \, s_i \, \varphi(s_i) \, ds_i\nonumber\\ 
=
     & \frac{1}{\sigma}\int_{-\infty}^{+\infty} \frac{\partial f(z_i, \bar x_i)}{\partial z_i}\, \frac{n z_i }{ds_i} \, \varphi(s_i) \, ds_i
    \nonumber\\ =
    &\int_{-\infty}^{+\infty} \frac{\partial f(z_i, \bar x_i)}{\partial z_i} \, \varphi(s_i) \, ds_i,
    \label{week}
\end{align}
\\
\noindent
where $z_i = x_i + \sigma s_i$.
By changing of variable, $s_i = \frac{z_i - x_i}{\sigma} $ we obtain that
\begin{equation}
    g_\sigma(x_i, \bar x_i) = \int_{-\infty}^{+\infty} \frac{\partial f(z_i, \bar x_i)}{\partial z_i} \,\frac{1}{\sigma}\varphi(\frac{z_i-x_i}{\sigma}) \, dz_i
\label{week2}
\end{equation}
and therefore, taking into account that a series of Gaussians $\frac{1}{\sigma_n}\varphi(\frac{z_i-x_i}{\sigma_n})$ with $\sigma_n\to 0$ defines a 
$\delta$-dirac distribution centered in $x_i$ \cite{gel2016generalized}, we have that
\begin{equation}
\lim_{\sigma\to 0} g_\sigma(x_i, \bar x_i) = \frac{\partial f(x)}{\partial x_i}.
\label{lim}
\end{equation}
\noindent

$\hfill \square$

Any entry of (\ref{qgrad})  is a weak definition of the derivative of $f(x)$ along $x_i$ \cite{gel2016generalized}.
Note that (\ref{week}) is well defined even though $f(x)$ is not differentiable at $(x_i,\,\bar x_i)$\footnote{Any $L_1$ function satisfying (\ref{week}), in place of $\frac{\partial f(z_i, \bar x_i)}{\partial z_i}$, is 
a \textit{weak derivative} of $f(x)$ along $x_i$.}.

\section{A New Estimate of the Gradient} \label{A new estimate of the gradient}
We consider the functional $g_\sigma(x_i, \bar x_i)$ which is the $i_{th}$ component of the gradient estimate (\ref{qgrad}) and, for the sake of simplicity, we write in a single formula the result of (\ref{week}) and (\ref{week2}).
\begin{align}
   g_\sigma(x_i, \bar x_i) = &\frac{1}{\sigma} \int_{-\infty}^{+\infty} 
f(x_i + \sigma s_i, \bar x_i) \, s_i \, \varphi(s_i) \, ds_i\nonumber\\
   = &\int_{-\infty}^{+\infty} \frac{\partial f(z_i, \bar x_i)}{\partial z_i} \,\frac{1}{\sigma}\varphi(\frac{z_i-x_i}{\sigma}) \, dz_i
    \label{week33}
\end{align}
Note that $\frac{1}{\sigma}\varphi(\frac{z_i-x_i}{\sigma})$ is $\mathcal{N}(x_i, \sigma^2)$.
Our goal consists in finding a numerical approximation of the first integral in (\ref{week33}). To do that, we compute the integral in a finite range, namely between \textit{-S} and \textit{S}
\begin{align}
 \Tilde{g}_\sigma(x_i, \bar x_i) := & \frac{1}{\sigma} \int_{-S}^{+S} f(x_i + \sigma s_i, \bar x_i) \, s_i \, \varphi(s_i) \, ds_i \nonumber\\
   = & - \frac{1}{\sigma} \int_{-S}^{+S} f(x_i + \sigma s_i, \bar x_i) \, \varphi'(s_i) \, ds_i
\label{truncated}
\end{align}
For $S$ sufficiently big the error between (\ref{week33}) and (\ref{truncated}) is negligible due to the fast decreasing of the Gaussian to infinity.
The definite integral in (\ref{truncated}) can be approximated by a quadrature formula, e.g. Trapezoidal Rule \cite{atkinson2008introduction}. Dividing the interval $[-S, S]$ in $2m$ sub-intervals, each of size $h = \frac{S}{m}$ we obtain:
\begin{align}
  \Tilde{g}_\sigma(x_i, \bar x_i) = -& \frac{h}{2\sigma} \bigg[ \bigg( f(x_i - \sigma S, \bar x_i)\,  \varphi'(-S) + f(x_i + \sigma S, \bar x_i) 
\, \varphi'(S)+ \nonumber   \\
  &2 \sum_{j = 1}^{2m-1} f(x_i + \sigma (-S + j h), \bar x_i)\, \varphi'(-S + j h) \bigg)\bigg] + \nonumber \\
 +& \frac{ h^2\,S}{6 \sigma} \frac{d}{d\tau^2} f(x_i + \sigma \tau, \bar x_i) \,  \varphi'(\tau) \hspace{15 pt} \tau \in [-S, S]. 
\label{cavalieri}
\end{align}

\noindent
It is well known that, under very general conditions, the trapezoidal quadrature formula (\ref{cavalieri})  has an error that is $\mathcal{O}(1/m^2)$ \cite{boyd2001chebyshev}. Indeed, once $\sigma$ and $S$ are chosen, we can easily check this property in our case. Let
\begin{align}
     \epsilon_\sigma (\tau,m) &= \frac{ h^2\,S}{6 \, \sigma} \frac{d}{d\tau^2} f(x_i + \sigma \tau, \bar x_i) \,  \varphi'(\tau) \hspace{15 pt} \tau \in [-S, S]. \\
     &= \frac{S^3}{6 \, \sigma\, m^2} \frac{d}{d\tau^2} f(x_i + \sigma \tau, \bar x_i) \,  \varphi'(\tau) \hspace{15 pt} \tau \in [-S, S]. \nonumber
\end{align}
    
\noindent Note that the derivatives of a guassian kernel $|\varphi^{(k)}(\tau)|$, up to the third order, are
all less than $1$ in absolute value for any $\tau$, and decrease rapidly as $\tau$ increases. Therefore, for $f$ sufficiently smooth in $(x_i\pm \sigma\,S)$, let
$$
K(x_i) = max\left( |f(x_i + \sigma \tau, \bar x_i)|,|\frac{d}{d\tau} f(x_i + \sigma \tau, \bar x_i)|,|\frac{d^2}{d\tau^2} f(x_i + \sigma \tau, \bar x_i)|\right).
$$
We can write:
$$
|\epsilon_\sigma (\tau,m)| \leq \frac{ h^2\,S}{6 \sigma} K(x_i) = \frac{ S^3 }{6 \sigma\, m^2} K(x_i),\quad \tau \in [-S, +S]
$$

\noindent Let us rewrite (\ref{cavalieri}) as follows
$$
\Tilde{g}_\sigma(x_i, \bar x_i) = \bar{g}_\sigma(x_i, \bar x_i) + \epsilon_\sigma (\tau,m)
$$
The larger the number of function evaluation $m$ the smaller the error term $\epsilon_\sigma (\tau,m)$.  On the other hand $\bar{g}_\sigma(x_i)$ can be interpreted as a combination of finite differences with some coefficients. Keeping in mind that $\varphi'(t) = -\varphi'(-t)$ and that $\varphi'(0) = 0$, after some simple algebra 
we can write:
\begin{align*}
  \bar{g}_\sigma(x_i, \bar x_i) = -\frac{h}{2\sigma}\bigg[|\varphi'(m\,h)| \bigg( f(x_i - \sigma\,m\,h, \bar x_i) - f(x_i + \sigma\,m\,h, \bar x_i)\bigg) +\\
  + 2 \sum_{j = 1}^{m-1} |\varphi'(jh)|  \bigg(f(x_i - \sigma jh, \bar x_i) - f(x_i + \sigma jh, \bar x_i)\bigg)\bigg]
\end{align*}
from which
\begin{align}
  \bar{g}_\sigma(x_i, \bar x_i) = \frac{h}{2\sigma}\bigg[|\varphi'(m\,h)|\, 2\sigma\,m\,h \,\frac{f(x_i + \sigma\,m\,h, \bar x_i) - f(x_i - \sigma\,m\,h, \bar x_i)}{2\sigma\,m\,h} +\nonumber\\
  +2\sum_{j = 1}^{m-1} |\varphi'(jh)|\,2\sigma\,j\,h \, \frac{f(x_i + \sigma jh, \bar x_i) - f(x_i - \sigma jh, \bar x_i)}{2\sigma\,j\,h} \bigg].
  \label{CCFD}
\end{align}
It is clear that $\bar{g}_\sigma(x_i, \bar x_i)$ is a linear combination of finite difference approximations, with different step sizes;
 for $\sigma h \to 0$ each one converges to the true value of the partial 
derivative ${\partial f(x_i, \bar x_i)}/{\partial x_i}$. Therefore the estimate $\bar{g}_\sigma(x_i, \bar x_i)$ converges to the true value only if the sum of its coefficients equals one. For this reason, it is advisable to \emph{normalize} the coefficients of the linear combination in (\ref{CCFD}) to eliminate the estimate bias for $\sigma$ finite. To this aim let $C$ be the sum of all the coefficients:
\begin{align}
&                       &&\negthickspace
\begin{rcases}
                C =  \sum_{j = 1}^{m} a^\prime_j, \\
       a^\prime_j = 2\,j\,h^2\,|\varphi^\prime(j h)|,\quad j=1,\ldots,m-1,\\
       a^\prime_m = m\,h^2\,|\varphi^\prime(m h)|,\\
\end{rcases}
       \label{aprime}
\end{align}

\medskip\noindent We can then write the normalized version of (\ref{CCFD}) as:
\begin{equation}
   \hat{g}_\sigma(x_i, \bar x_i) = \sum_{j = 1}^{m} a_j\, \frac{ f(x_i + \sigma\,j\,h, \bar x_i) - f(x_i - \sigma\,j\,h, \bar x_i)} {2\sigma\,j\,h}
  \label{CCFD_norm}
\end{equation}
where
\begin{equation}
    a_j = \, \frac{a^\prime_j}{C},\qquad \sum_{j = 1}^{m} a_j = 1.
    \label{cstnt}
\end{equation}

\noindent For $\sigma$ small enough the normalization of the coefficients may not be necessary, the distorsion of the estimate being negligible. Let us now evaluate the error bound corresponding to estimate (\ref{CCFD_norm}), from here on referred to as NMXFD (Normalized Mixed Finite Difference).
\begin{theorem}
\label{error_bounds}
Let $f(x)$ be twice continuously differentiable  and its Hessian be $H$-Lipschitz for all $x\in R^n$. Consider the gradient approximation obtained 
by (\ref{CCFD_norm})
\begin{equation}
    \widehat G_\sigma (x) = \left[\, \hat{g}_\sigma(x_1),\ldots \hat{g}_\sigma(x_n)\,\right]^T
    \label{NMXFD_vettore}
\end{equation}

We have that
\begin{equation}
    \|\widehat G_\sigma(x) - \nabla f(x)\| \leq \sqrt{n}\>\frac{H \sigma^2\,S^2 }{6}
    \label{bounds_norm}
\end{equation}
\end{theorem}
\textbf{Proof}: Any single finite difference term in (\ref{CCFD_norm}) has an error with respect to the true value ${\partial f(x_i, \bar x_i)}/{\partial x_i}$ whose bound depends on the step size and on the regularity properties of function $f$. From \cite{berahas2019theoretical} we have that
\begin{equation}
    \left |\frac{ f(x_i + \sigma\,j\,h, \bar x_i) - f(x_i - \sigma\,j\,h, \bar x_i)} {2\sigma\,j\,h} - \frac{\partial f(x_i, \bar x_i)}{\partial x_i}\right |\leq \frac{H \sigma^2 (jh)^2}{6}
    \label{cfderror}
\end{equation}

\noindent for $j=1,\ldots,m$. Therefore, since $\sum_{j = 1}^{m} a_j =  1$, and $a_j>0$, $j=1,\ldots,m$, we can write
\begin{align}
    &\left |\hat{g}_\sigma(x_i) - \frac{\partial f(x_i, \bar x_i)}{\partial x_i}\right | = \left |\hat{g}_\sigma(x_i) - \sum_{j=1}^m a_j\,\frac{\partial f(x_i, \bar x_i)}{\partial x_i}\right |
    \leq\nonumber\\
    & \sum_{j = 1}^{m} a_j\, \left | \frac{ f(x_i + \sigma\,j\,h, \bar x_i) - f(x_i - \sigma\,j\,h, \bar x_i)} {2\sigma\,j\,h}-\frac{\partial f(x_i, \bar x_i)}{\partial x_i} \right |\leq\nonumber\\
    &\frac{H \sigma^2\,h^2 }{6}\left(\sum_{j=1}^{m} a_j \,j^2\right)\leq \frac{H \sigma^2\,h^2\,m^2 }{6} = \frac{H \sigma^2\,S^2 }{6}. \nonumber
\end{align}
which, applied to all entries of $\widehat G_\sigma(x)-\nabla f(x)$, proves the theorem. $\hfill \square $
Here we used the equality $m\,h=S$ that implies that the error bound does not depend on the number of function evaluations.

\section{Estimation Error with Noisy Data}
\label{Noisy_data}
Let us now evaluate how the performance of the gradient estimate $NMXFD$ (\ref{NMXFD_vettore}) here referred to as $\hat{G}_\sigma^{\text{\tiny{MXF}}}(x)$ compares with that of the Central Finite Differences ($CFD$), taking also into account the presence of an additive noise affecting  the sampled function values $f(x)$. 
Let $\{e_i\}$ be the canonical base of $R^n$, then we can write:
\begin{equation}
\hat{G}_\sigma^{\text{\tiny{MXF}}}(x) = \sum_{i=1}^n \hat{g}_\sigma(x_i)\, e_i
\label{EMXF}
\end{equation}
\noindent
With the same notation we can easily write the gradient estimate according to the CFD scheme here denoted as $\hat{G}_\sigma^{\text{\tiny{CFD}}}(x)$:
\begin{equation}
\hat{G}_\sigma^{\text{\tiny{CFD}}}(x) = \sum_{i=1}^n \frac{ f(x_i + \sigma\,h, \bar x_i) - f(x_i - \sigma\,h, \bar x_i)}{2\sigma\,h}\,e_i = \sum_{i=1}^n \delta f_\sigma(x_i)\, e_i
\label{ECFD}
\end{equation}

\noindent Let $\{\epsilon_i\}$ denote a discrete random field modeling the additive noise on the sampled function values with the following properties: $\epsilon_i \sim N(0,\lambda^2)$ and $E[\epsilon_i \,\epsilon_j] = 0$ for $i\ne j$.
\noindent
We now compute the estimation errors for the two schemes and compare them in 
terms of accuracy (mean value) and precision (variance). The accuracy evaluates the estimate bias, i.e. the systematic source of the error, like the limited the number $N$ of function evaluations used to build the estimate. The precision is the dispersion of the estimation error around its mean value and evaluates the variability of the statistic source of the error.

\newpage
\noindent\textit{The CFD scheme}

\noindent
According to (\ref{ECFD}), a number $N = 2 n$ of function evaluations is considered to obtain
$$
\hat{G}_\sigma^{\text{\tiny{CFD}}}(x) = \sum_{i=1}^n \delta f_\sigma(x_i)\, e_i +\sum_{i=1}^n 
\frac{\epsilon_i^+ - \epsilon_i^-}{2\,\sigma\,h} \, e_i
$$
with $\epsilon_i^\pm$ denoting the noise on the function values  $f_\sigma(x_i\pm\sigma\,h,\bar x_i)$. Let 
$$e_{\text{\tiny{CFD}}}(x) = \hat{G}_\sigma^{\text{\tiny{CFD}}}(x) - \nabla f(x)$$ 
be the estimation error. We can see that

$$
E[e_{\text{\tiny{CFD}}}(x)] =  \sum_{i=1}^n \delta f_\sigma(x_i)\, e_i
$$
and
\begin{equation}
\label{varCFD}
    var[e_{\text{\tiny{CFD}}}(x)] =  n\, \frac{2 \lambda^2}{4\, \sigma^2 \,h^2} = \frac{n\,\lambda^2}{2\sigma^2 \,h^2}
\end{equation}

\noindent where $var[z]$, $z \in R^n$ with $E[z] = 0$, indicates the trace of the 
covariance matrix $E[z\,z^T]$.
\noindent
Now, for functions $f$ as in theorem (\ref{error_bounds}), let us consider the property (\ref{cfderror}), with $j=1$, for all the components of $E[e_{\text{\tiny{CFD}}}(x)]$. We obtain that
$$
\|E[e_{\text{\tiny{CFD}}}(x)]\|\le \sqrt n\,\frac{H\,\sigma^2\,h^2}{6}.
$$
Therefore as the increment $\sigma h\to 0$, the error goes to zero as well \textit{on average}, but its variance increases without bound as $\mathcal{O}\left(1/(\sigma h)^2\right)$.
\vskip 1cm
\noindent\textit{The NMXFD scheme}\\
\noindent
In this case, according to (\ref{EMXF}), a number $N = 2 m\,n$ of function evaluations is considered 
to obtain 
$$
\hat{G}_\sigma^{\text{\tiny{MXF}}}(x) = \sum_{i=1}^n \hat{g}_\sigma(x_i)\, e_i + \sum_{i=1}^n\left( \sum_{j = 1}^{m} a_j\, \frac{\epsilon_{i,j}^+ - \epsilon_{i,j}^-} {2\sigma\,j\,h}\right)\,e_i
$$
with $\epsilon_{i,j}^\pm$ denoting the error terms on the function values $f(x_i\pm\sigma\,jh, \bar x_i)$, $i=1,\ldots,n$, $j=1,\ldots,m$. For the estimation error 
$$e_{\text{\tiny{MXF}}}(x) = \hat{G}_\sigma^{\text{\tiny{MXF}}}(x) - \nabla f(x),$$ 
we readily obtain that
\begin{align}
    E\left[e_{\text{\tiny{MXF}}}(x)\right] &= \sum_{i=1}^n \hat{g}_\sigma(x_i)\, e_i-\nabla f(x) \nonumber\\
    var\left[e_{\text{\tiny{MXF}}}(x)\right] &= \frac{n\,\lambda^2}{2\sigma^2h^2}\left(\sum_{j=1}^{m} \frac{a_j^2}{j^2}\right) \label{varMXF}
\end{align}

\noindent Under the assumptions of theorem (\ref{error_bounds}), and taking into account (\ref{bounds_norm}), we obtain
\begin{equation}
    \| E\left[e_{\text{\tiny{MXF}}}(x)\right] \| \leq \sqrt{n}\>\frac{H \sigma^2\,m^2\,h^2 }{6}.
    \label{biasMXF}
\end{equation}

\noindent As for the error variance two interesting results can be proved. 

\begin{proposition} For any $m>1$, the variance of the estimation error of the $NMXFD$ scheme is strictly lower than the variance of the estimation error of the $CFD$ scheme, i.e.
\begin{equation}
    var\left[e_{\text{\tiny{MXF}}}(x)\right] < var\left[e_{\text{\tiny{CFD}}}(x)\right]
\end{equation}
in any $x\in R^n$ and for any $\sigma$, $h$.
\end{proposition}

\noindent\textbf{Proof.}
The sum of squares $\sum_{j= 1}^{m} a_j^2$ is strictly less then 1 since the coefficients $a_j$, $j=1,\ldots,m$, are all positive and their sum is 1.

Therefore from (\ref{varMXF}) we obtain that 
\begin{equation}
    var\left[e_{\text{\tiny{MXF}}}(x)\right] = \frac{n \lambda^2}{2 \sigma^2\,h^2}  \sum_{j = 1}^{m} \frac{a_j^2}{j^2} < \frac{n \lambda^2}{2 \sigma^2\,h^2} =  var\left[e_{\text{\tiny{CFD}}}(x)\right] .
\end{equation}
$\hfill \square$

\noindent
Now we further show that $var\left[e_{\text{\tiny{MXF}}}(x)\right]$ goes to zero as $N$ increases. 
\begin{proposition} For any $x\in R^n$, the variance of the estimation error of the $NMXFD$ scheme has the following asymptotic behavior
\begin{equation}
    var\left[e_{\text{\tiny{MXF}}}(x)\right] \sim  \mathcal{O}\left(\frac{1}{N}\right).
\end{equation}
\end{proposition}
\noindent{\textbf Proof.} By taking into account relations (\ref{aprime}),  we have that
\begin{align}
    C =& m\,h^2\left(|\varphi^\prime(m h)|+2\sum_{j=1}^{m-1}\frac{j}{m}|\varphi^\prime(j h)|\right)\nonumber\\
    &\le 2m\,h\,\frac{h}{2}\left(|\varphi^\prime(m h)|+2\sum_{j=1}^{m-1}|\varphi^\prime(j h)|\right).
    \label{C}
\end{align}
Let us denote with $I_{\varphi^\prime}^{(1)}(m)$ the following quantity
$$
I_{\varphi^\prime}^{(1)}(m) = \frac{h}{2}\left(|\varphi^\prime(m h)|+2\sum_{j=1}^{m-1}|\varphi^\prime(j h)|\right)
$$
that is the trapezoidal quadrature formula for the integral
$$
\int_0^S |\varphi^\prime(t)|\,dt  = \frac{1}{\sqrt{2\pi}}\left(1-e^{-\frac{S^2}{2}}\right).
$$
Due to the $\mathcal{O}(1/N^2) $  property of the error of the trapezoidal 
rule, we have that
$$ \bigg|I_{\varphi^\prime}^{(1)}(m) - \frac{1}{\sqrt{2\pi}}\left(1-e^{-\frac{S^2}{2}}\right)\bigg|= \mathcal{O}(1/N^2).
$$
Therefore, from (\ref{C}), we easily obtain that
\begin{equation}
\bigg| C-\frac{2m\,h}{\sqrt{2\pi}}\left(1-e^{-\frac{S^2}{2}}\right)\bigg|\le 2m\,h\,\bigg|I_{\varphi^\prime}^{(1)}(m) - \frac{1}{\sqrt{2\pi}}\left(1-e^{-\frac{S^2}{2}}\right)\bigg| = \mathcal{O}(1/N^2) 
\label{CP}
\end{equation}
so that $C$ is a bounded quantity as $N=2m\,n$ increases (by increasing $m$), taking into account that $mh=S$.
\noindent Now, according to the relations (\ref{cstnt}) we can write
\begin{align*}
\sum_{j=1}^m \frac{a_j^2}{j^2} &= \frac{1}{C^2}\left( \frac{m^2\,h^4|\varphi^\prime(m h)|^2}{m^2}+\sum_{j=1}^{m-1}\frac{4\,j^2\,h^4\,|\varphi^\prime(j h)|^2}{j^2} \right) \nonumber \\
&= \frac{h^4}{C^2}\left(|\varphi^\prime(m h)|^2 +2\sum_{j=1}^{m-1} 2|\varphi^\prime(j h)|^2 \right)\nonumber\\
& \leq \frac{2\,h^3}{C^2}\frac{h}{2}\left(2|\varphi^\prime(m h)|^2 +2\sum_{j=1}^{m-1} 2|\varphi^\prime(j h)|^2 \right).\nonumber\\
\end{align*}
Define now $I_{\varphi^\prime}^{(2)}(m)$ as follows
$$
I_{\varphi^\prime}^{(2)}(m) = \frac{h}{2}\left(2|\varphi^\prime(m h)|^2 
+2\sum_{j=1}^{m-1} 2|\varphi^\prime(j h)|^2 \right).
$$
It is  the trapezoidal quadrature rule for the integral
$$
2\int_0^S |\varphi^\prime(t)|^2\,dt = \sqrt\pi\, \erf(S) - S e^{-S^2} = \Phi(S),
$$
where $\erf(z)=\frac{2}{\sqrt\pi}\int_0^z e^{-t^2}\,dt$ is the Gauss error function. Hence, for the usual property of the error, we can write
$$
\bigg|I_{\varphi^\prime}^{(2)}(m) - \Phi(S)\bigg| =  \mathcal{O}(1/N^2).
$$
Therefore we obtain that
\begin{align*}
  var\left[e_{\text{\tiny{MXF}}}(x)\right] =& \frac{n\,\lambda^2}{2\sigma^2\,h^2}\left(\sum_{j=1}^{m} \frac{a_j^2}{j^2}\right)\leq \frac{n\,\lambda^2}{2\sigma^2\,h^2} \frac{2\,h^3}{C^2}I_{\varphi^\prime}^{(2)}(m) \\
  \leq &\frac{n\,\lambda^2}{\sigma^2}\,\frac{h}{C^2}\left(\big|I_{\varphi^\prime}^{(2)}(m)- \Phi(S)\big| + \big|\Phi(S) \big|\right). \\
\end{align*}

\noindent Now recalling that $m\,h = S $, and that $N=2m\,n$, we can write
\begin{align*}
 var\left[e_{\text{\tiny{MXF}}}(x)\right]\leq & \frac{n\,\lambda^2}{\sigma^2}\,\frac{S}{m\,C^2}\left(\big|I_{\varphi^\prime}^{(2)}(m)- \Phi(S)\big| + \big|\Phi(S) \big|\right)\\
 \leq & \frac{2}{N}\,\frac{n^2\,\lambda^2\,S}{\sigma^2\,C^2}\left(\big|I_{\varphi^\prime}^{(2)}(m)- \Phi(S)\big| + \big|\Phi(S) \big|\right)\\
\end{align*}
which, along with (\ref{CP}), proves the proposition.
$\hfill \square$


\noindent

\section{Numerical Experiments}
\label{Numerical experiment}
We tested our method for estimating the gradient by comparing its performance with those of other methods on 69 functions from the Schittkowski test set \cite{schittkowski2012more}.

\noindent
For each function we did the following: we generated a random starting point $x^0$, and minimized the function using the quasi-Newton method of Broyden, Fletcher, Goldfarb, and Shanno (BFGS) \cite{nocedal2006sequential}, finding the optimal point $x^*$ with $\nabla f(x^*) \approx 0$. We then identified the first instance of a point $x^k$ where
$$
\frac{\|\nabla f(x^k)\|}{\|\nabla f(x^0)\|} \leq \alpha
$$
for each of the following values of $\alpha$: $10^{0},
10^{-1},10^{-2},10^{-3},10^{-4},10^{-5},10^{-6}$.
\noindent
In this way we generated seven different buckets, one for each $\alpha$, of 69 different points, one for each function. Bucket \textit{i} indicates the one associated to $\alpha = 10^{-i}$. Bucket 0 is therefore the one with the points that are further from the optimal solution and bucket 6 is the one with points closer to the optimal solution.

\noindent
Then, for each point we computed the gradient approximations obtained with the Normalized MiXed Finite Differences scheme (NMXFD) and with those considered benchmarks by the literature, namely: Forward Finite Differences (FFD), Central Finite Differences (CFD), Gaussian Smoothed Gradient (GSG), Central Gaussian Gmoothed Gradient (cGSG) as defined in \cite{berahas2019theoretical}.
\noindent
Different tables will summarize the results of this comparison.\\
The tables show, for different values of the number of function evaluations ($N$) and different buckets ($B$), the median value of the log of the relative approximation error over all the 69 points in each bucket. 
 
\noindent
We define relative approximation error as $$
\eta = \frac{\|g(x)- \nabla f(x)\|}{\|\nabla f(x)\|}
$$ where $g(x)$ is the generic gradient estimate.
The number of function evaluations $N$ is expressed in the following tables as a function of the number of dimensions $n$.
$FFD$ and $CFD$ schemes only allow for a specific value of \textit{N} ($n+1$ and $2n$, respectively).
In $GSG$ and in $cGSG$, N is linked to the number of direction sampled to build the gradient approximation ($N=(M+1)$ in (\ref{fd}) and $N = 2\,M$ in (\ref{cfd})). In the NMFXD scheme, the value of $N$ is linked to the value of $m$ in formula (\ref{CCFD_norm}). In particular, we have that $N = 2m\,n$.
\noindent
In each table, the lowest entry for every bucket is highlighted in bold, and the second lowest is underlined. 

\noindent
\subsection{Noise-free setting}
For the noise-free setting we report three different tables obtained using a different value of $\sigma$ (shared by all the schemes) to compute the gradient approximation.

\begin{table}[htbp]
  \centering
    \renewcommand{\arraystretch}{0.70} 
  \caption{Median log of relative error with \textbf{$\sigma = 10^{-2}$}}
    \begin{tabular}{ccrrrrrrr}
    Scheme & N     & \multicolumn{1}{c}{B0} & \multicolumn{1}{c}{B1} & \multicolumn{1}{c}{B2} & \multicolumn{1}{c}{B3} & \multicolumn{1}{c}{B4} & \multicolumn{1}{c}{B5} & \multicolumn{1}{c}{B6} \\
    \midrule
    \textbf{FFD} & n+1   & 0.08  & 1.22  & 2.20  & 3.43  & 4.43  & 5.47  & 6.52 \\
    \midrule
    \textbf{CFD} & 2n    & \textbf{-2.26} & \textbf{-1.13} & \textbf{-0.32} & \textbf{0.69 } & \textbf{1.60}  & \textbf{2.41}  &\textbf{ 3.58} \\
    \midrule
          & 2n+1  & 1.84  & 1.92  & 2.52  & 3.57  & 4.69  & 5.63  & 6.92 \\
    \textbf{GSG} & 4n+1  & 1.70  & 1.78  & 2.34  & 3.41  & 4.50  & 5.46  & 6.82 \\
          & 8n+1  & 1.54  & 1.66  & 2.10  & 3.31  & 4.49  & 5.38  & 6.67 \\
    \midrule
          & 2n    & 1.97  & 1.96  & 1.96  & 1.99  & 2.08  & 2.44  & 3.46 \\
    \textbf{cGSG} & 4n    & 1.86  & 1.82  & 1.86  & 1.90  & 2.06  & 2.95  & 4.01 \\
          & 8n    & 1.71  & 1.69  & 1.74  & 1.81  & 2.13  & 3.14  & 4.28 \\
    \midrule
          & 2n    & -1.66 & -0.53 & 0.28  & 1.29  & 2.26  & 3.06  & 4.18 \\
    \textbf{NMXFD} & 4n    & -1.98 & -0.84 & -0.03 & 0.97  & 1.92  & 2.71  & 3.87 \\
          & 8n    & \underline{-1.99} & \underline{-0.86} & \underline{-0.05} & \underline{0.96}  & \underline{1.90}  & \underline{2.68}  & \underline{3.85} \\
    \bottomrule
    \end{tabular}%
  \label{tab:no_noise_2}%
\end{table}%

\begin{table}[htbp]
  \centering
    \renewcommand{\arraystretch}{0.70} 
  \caption{Median log of relative error with \textbf{$\sigma = 10^{-5}$}}
    \begin{tabular}{ccrrrrrrr}
    Scheme & N     & \multicolumn{1}{c}{B0} & \multicolumn{1}{c}{B1} & \multicolumn{1}{c}{B2} & \multicolumn{1}{c}{B3} & \multicolumn{1}{c}{B4} & \multicolumn{1}{c}{B5} & \multicolumn{1}{c}{B6} \\
    \midrule
    \textbf{FFD} & n+1   & -2.92 & -1.78 & -0.80 & 0.43  & 1.43  & 2.47  & 3.51 \\
    \midrule
    \textbf{CFD} & 2n    & \textbf{-8.17} & \textbf{-6.99} & \textbf{-6.14} & \textbf{-4.87} & \textbf{-3.75} & \textbf{-3.07} & \textbf{-1.73} \\
    \midrule
          & 2n+1  & 1.84  & 1.84  & 1.85  & 1.84  & 2.00  & 2.64  & 3.92 \\
    \textbf{GSG} & 4n+1  & 1.69  & 1.71  & 1.71  & 1.74  & 1.90  & 2.48  & 3.80 \\
          & 8n+1  & 1.53  & 1.57  & 1.56  & 1.57  & 1.77  & 2.41  & 3.67 \\
    \midrule
          & 2n    & 1.96  & 1.96  & 1.96  & 1.96  & 1.96  & 1.97  & 1.93 \\
    \textbf{cGSG} & 4n    & 1.86  & 1.82  & 1.85  & 1.85  & 1.85  & 1.85  & 1.83 \\
          & 8n    & 1.71  & 1.68  & 1.70  & 1.68  & 1.71  & 1.71  & 1.71 \\
    \midrule
          & 2n    & -7.66 & -6.47 & -5.58 & -4.43 & -3.24 & -2.75 & -1.21 \\
    \textbf{NMXFD} & 4n    & -7.90 & -6.74 & -5.85 & -4.67 & -3.55 & -2.79 & -1.45 \\
          & 8n    & \underline{-7.95} & \underline{-6.76 } & \underline{-5.87} & \underline{-4.73} & \underline{-3.57} & \underline{-2.84} & \underline{-1.54} \\
    \bottomrule
    \end{tabular}%
  \label{tab:no_noise_5}%
\end{table}%

\begin{table}[htbp]
  \centering
    \renewcommand{\arraystretch}{0.70} 
  \caption{Median log of relative error with \textbf{$\sigma = 10^{-8}$}}
    \begin{tabular}{ccrrrrrrr}
    Scheme & N     & \multicolumn{1}{c}{B0} & \multicolumn{1}{c}{B1} & \multicolumn{1}{c}{B2} & \multicolumn{1}{c}{B3} & \multicolumn{1}{c}{B4} & \multicolumn{1}{c}{B5} & \multicolumn{1}{c}{B6} \\
    \midrule
    \textbf{FFD} & n+1   & -5.56 & -4.73 & -3.74 & 1.43  & -1.50 & -0.44 & 0.67 \\
    \midrule
    \textbf{CFD} & 2n    & -6.00 & -6.20 & -6.23 &\textbf{ -3.75} & -6.20 & -6.25 & \textbf{-6.23} \\
    \midrule
          & 2n+1  & 1.84  & 1.84  & 1.84  & 2.00  & 1.82  & 1.83  & 1.91 \\
    \textbf{GSG} & 4n+1  & 1.69  & 1.71  & 1.72  & 1.90  & 1.70  & 1.69  & 1.79 \\
          & 8n+1  & 1.53  & 1.57  & 1.56  & 1.77  & 1.55  & 1.55  & 1.66 \\
    \midrule
          & 2n    & 1.96  & 1.96  & 1.96  & 1.96  & 1.96  & 1.96  & 1.93 \\
    \textbf{cGSG} & 4n    & 1.86  & 1.82  & 1.85  & 1.85  & 1.84  & 1.85  & 1.82 \\
          & 8n    & 1.71  & 1.68  & 1.70  & 1.71  & 1.71  & 1.71  & 1.71 \\
    \midrule
          & 2n    & -6.48 & -6.36 & -6.52 & -3.24 & -6.41 & -6.42 & -6.09 \\
    \textbf{NMXFD} & 4n    & \underline{-6.17} & \underline{-6.29} & \underline{-6.41} & -3.55 & \underline{-6.43} & \underline{-6.48} & \underline{-6.20} \\
          & 8n    & \textbf{-6.42} &\textbf{ -6.44} & \textbf{-6.44} & \underline{-3.57} & \textbf{-6.50} & \textbf{-6.51} & -6.15 \\
    \bottomrule
    \end{tabular}%
  \label{tab:no_noise_8}%
\end{table}%

\noindent
It is possible to notice that in a noise-free setting, lower values of $\sigma$ tend to yield to better results, as one would expect from the theory. 
The closer the point is to the minimum value of a function, the harder it is to obtain an accurate estimate of its gradient, unless $\sigma$ is very small. As a matter of fact, for points belonging to lower index buckets - thus far from the minimum of the function, the value $\sigma = 10^{-5}$ yields the better performances, while accurate estimates of the gradient of points closer to the minimum value of a function require using of a lower value of $\sigma$.
We can also see that the error of the proposed method, NMXFD, is of the same order of magnitude of that of CFD, and almost always better than that of the other methods.

In our experiments, we have also produced gradient estimates using two more methods: 
\begin{itemize}
    \item[$\bullet$] by removing the normalization of the coefficients in the computation of NMXFD, i.e. implementing the gradient approximation as in (\ref{CCFD}).  
    \item[$\bullet$] by computing the estimate as the raw average of central finite differences at different stepsizes, that is (\ref{CCFD_norm}) with $a_j = \frac{1}{m}$.
\end{itemize}
Both of these methods performed consistently worse than NMXFD, and they have not been reported in the tables for brevity.
Still, the better performances of NMXFD over the raw average of central finite differences seem to confirm that the rationale behind the choice of coefficients used to weight the CFDs in the proposed approach is promising from a computational point of view.

\subsection{Noisy setting}
\noindent
We also show results of the noisy scenario, where the noise term is described in Section 4 and has $\lambda = 0.001$.
\noindent
The estimation procedure is slightly different from the one of the noise-free setting. In table \ref{tab:noise_2}, the median log of the relative errors $\eta_i$ of the 69 different Schittkowski function is reported. Each $\eta_i$ is computed as the average of 100 relative approximation errors, resulting from 100 independent noise realizations. The rationale behind this choice was to mitigate the dependence of the results from one particular noise realization.
\noindent
Results are shown in table \ref{tab:noise_2}, where the gradient estimates are obtained with $\sigma = 0.01$. 

\begin{table}[htbp]
  \centering
    \renewcommand{\arraystretch}{0.70} 
  \caption{Median log of relative error with \textbf{$\sigma = 10^{-2}$}, noisy setting.}

    \begin{tabular}{ccrrrrrrr}
    Scheme & N     & \multicolumn{1}{c}{B0} & \multicolumn{1}{c}{B1} & \multicolumn{1}{c}{B2} & \multicolumn{1}{c}{B3} & \multicolumn{1}{c}{B4} & \multicolumn{1}{c}{B5} & \multicolumn{1}{c}{B6} \\
    \midrule
    \textbf{FFD} & n+1   & 0.22  & 1.45  & 2.46  & 3.57  & 4.86  & 5.74  & 6.86 \\
    \midrule
    \textbf{CFD} & 2n    & -1.06 & 0.10  & 1.34  & 2.23  & 3.50  & 4.32  & 5.49 \\
    \midrule
          & 4n+1  & 1.72  & 1.79  & 2.56  & 3.66  & 4.82  & 5.69  & 6.84 \\
    \textbf{GSG} & 8n+1  & 1.56  & 1.66  & 2.43  & 3.51  & 4.67  & 5.56  & 6.70 \\
          & 12n+1 & 1.47  & 1.56  & 2.35  & 3.43  & 4.59  & 5.48  & 6.61 \\
    \midrule
          & 4n    & 1.85  & 1.85  & 1.86  & 2.39  & 3.61  & 4.33  & 5.65 \\
    \textbf{cGSG} & 8n    & 1.71  & 1.71  & 1.73  & 2.29  & 3.55  & 4.28  & 5.61 \\
          & 12n   & 1.62  & 1.63  & 1.65  & 2.24  & 3.52  & 4.25  & 5.58 \\
    \midrule
          & 4n    & -1.20 & 0.00  & 1.17  & 2.23  & 3.43  & 4.24  & 5.52 \\
    \textbf{NMXFD} & 8n    & \underline{-1.31}  & \underline{-0.09}  & \underline{1.08}   & \underline{2.15}   & \underline{3.40}   & \underline{4.19}   & \underline{5.42}  \\
          & 12n   & \textbf{-1.36} & \textbf{-0.15} & \textbf{1.05}  & \textbf{2.11}  & \textbf{3.39}  & \textbf{4.15}  &\textbf{ 5.38} \\
    \bottomrule

    \end{tabular}%
  \label{tab:noise_2}%
\end{table}%

\noindent
Table \ref{tab:noise_2} shows that NMXFD performs better than the other schemes in presence of noise, although reasonably low relative approximation errors are obtained only for the first three buckets. For the other ones the error $\eta$ increases significantly. This is due to the fact that the denominator of $\eta$ gets smaller as we move to points close to the minimum value of the function, while the variance of the approximation error does not change across different buckets.
Just like in the noise-free setting, increasing the number of function evaluations allows to increase the precision of all the schemes, as expected from the theory.

\noindent
Different values of $\sigma$ for estimating the gradient ($10^{-1}$ , $10^{-3}$, $10^{-4}$) have also been used. The associated tables have not been reported for brevity, since they yielded to the same conclusions and since the performances for almost every method and every bucket with those values of $\sigma$ are significantly worse. This can be inferred from the theory, since the value of $\sigma$ influences the bias and the variance of the estimate error in opposite directions, as we can see from (\ref{varMXF}) and (\ref{biasMXF})  in section \ref{Noisy_data}.

\noindent
Data availability statement: Data sharing not applicable to this article as no datasets were generated or analysed during the current study.

\section{Conclusions}
\label{conclusioni} In this paper a novel scheme to estimate the gradient of a function is proposed. It is based on linear functionals defining a filtered version of the objective function. Unlike standard methods where the approximation error is characterized from a statistical point of view and therefore may be quite large on a given experiment, one advantage of the proposed scheme relies on a deterministic characterization of the approximation error in the noise-free setting.

\noindent
The other advantage lies in its behaviour when function evaluations are affected by noise. In fact, the variance of the estimation error of the proposed method is showed to be strictly lower than that of the Central Finite Difference scheme, and diminishes as the number of function evaluations increases. The suitable linear combination of finite differences seems to have a filtering role in the case of noisy functions, thus resulting in a more robust estimator.

Numerical experiments on a significant benchmark given by the 69 Schittkowski functions show the good performances of the proposed method when compared with those of the standard methods commonly used in the literature. In particular, the performances of NMXFD are comparable with those of CFD in absence of noise and better with noisy data, and seem to be better than those of other schemes in both scenarios. The impact of the use of the proposed gradient approximation in optimization algorithms has not been studied yet and will be developed in a future study.

\appendix  
\section*{Appendix}
\textbf{Proof of theorem (\ref{prova_teorema0}): \\
}
we have that
$$
|| G_{\sigma}(x) - \overline G_\sigma(x) ||^2=\sum_{i=1}^n \left((G_{\sigma}(x))_i - (\overline G_\sigma(x))_i\right)^2
$$
where $(G_{\sigma}(x))_i$ is given  by (\ref{f6})
$$
(G_{\sigma}(x))_i = \int_{R^{n-1}} g_\sigma(x_i,\,\bar x_i + \sigma \bar s_i)\,\varphi(\bar s_i)\,d\bar s_i
$$
and $(\overline G_\sigma(x))_i = g_\sigma(x_i,\,\bar x_i )$, by (\ref{single_component}). We can write
\begin{align}
    \left((G_{\sigma}(x))_i - (\overline G_\sigma(x))_i\right)^2 =& \left( \int_{R^{n-1}} g_\sigma(x_i,\,\bar x_i + \sigma \bar s_i)\,\varphi(\bar s_i)\,d\bar s_i - g_\sigma(x_i,\,\bar x_i )\right)^2 \nonumber \\
    =&\left( \int_{R^{n-1}} (g_\sigma(x_i,\,\bar x_i + \sigma \bar s_i)-g_\sigma(x_i,\,\bar x_i ))\,\varphi(\bar s_i)\,d\bar s_i \right)^2 \label{niff_tra_g}
\end{align}
where the last equality holds since $\int_{R^{n-1}} \varphi(\bar s_i)\,d\bar s_i = 1$.

Now, the integrand in (\ref{niff_tra_g}) has the following expression
\begin{align}
    &g_\sigma(x_i,\,\bar x_i + \sigma \bar s_i)-g_\sigma(x_i,\,\bar x_i ) 
=\nonumber \\
    &\frac{1}{\sigma}\int_{-\infty}^\infty \left( f(x_i + \sigma s_i,\,\bar x_i + \sigma \bar s_i)-f(x_i + \sigma s_i,\,\bar x_i)\right)\,s_i\,\varphi(s_i) \, d s_i,
    \label{square_error}
\end{align}
and for the argument of the integral we can write

\begin{align}
&f(x_i + \sigma s_i,\,\bar x_i + \sigma \bar s_i)-f(x_i + \sigma s_i,\,\bar x_i)= \nonumber \\
=& \left(f(x_i + \sigma s_i,\,\bar x_i + \sigma \bar s_i)-f(x_i,\,\bar x_i)\right)-\left(f(x_i + \sigma s_i,\,\bar x_i)-f(x_i,\,\bar x_i)\right)\nonumber \\
=& \,\nabla f(x^\prime)^T\,\sigma s -(\nabla f(x_i^{\prime\prime}, \bar{x_i}))_i\,\sigma\,s_i \nonumber\\
=&\,(\nabla f(x^\prime))_i\,\sigma\,s_i+\overline{(\nabla f(x^\prime))}_i^T\,\sigma\,\bar s_i-(\nabla f(x_i^{\prime\prime}, \bar{x_i}))_i\,\sigma\,s_i \label{gdiff}
\end{align}
with $x^\prime\in(x,x+\sigma s)$ and $x_i^{\prime\prime} \in(x_i,x_i+\sigma s_i)$. \\
We further have that
\begin{equation}
\overline{(\nabla f(x^\prime))}_i = \overline{(\nabla f(x^\prime))}_i - 
\overline{(\nabla f(x))}_i + \overline{(\nabla f(x))}_i
\label{nabla_barrato}
\end{equation}

Now substituting (\ref{gdiff}) and (\ref{nabla_barrato}) into (\ref{square_error}) we obtain that
\begin{align}
  &g_\sigma(x_i,\,\bar x_i + \sigma \bar s_i)-g_\sigma(x_i,\,\bar x_i ) =\nonumber\\
  &=  \int_{-\infty}^\infty \left[ (\nabla f(x^\prime))_i - (\nabla f(x_i^{\prime\prime}, \bar{x_i}))_i\right]\,s_i^2\,\varphi(s_i) \, d s_i  + \nonumber\\
 &+ \int_{-\infty}^\infty \left(\overline{(\nabla f(x^\prime))}_i -\overline{(\nabla f(x))}_i + \overline{(\nabla f(x))}_i   \right)^T\bar{s_i}\,s_i\,\varphi(s_i) \, d s_i.\label{quasi}
\end{align}
By the Lipschitz property of the gradient, and recalling that
$$
\int_{-\infty}^\infty \overline{(\nabla f(x))}_i^T\,\bar{s_i}\,s_i\,\varphi(s_i) \, d s_i = 0
$$
we have:
\begin{align}
 & \left(g_\sigma(x_i,\,\bar x_i + \sigma \bar s_i)-g_\sigma(x_i,\,\bar x_i )\right)^2 \leq \nonumber \\
& \int_{-\infty}^{\infty} L^2 \sigma^2 \, \|s\|^2 \, s_i^4\,\varphi(s_i)\, d s_i +  \nonumber \\
& \int_{-\infty}^{\infty} L^2 \sigma^2 \, \|s\|^2 \,\|\bar{s_i}\|^2 \, s_i^2\,\varphi(s_i)\, d s_i
\label{qquasi}
\end{align}
We can finally substitute (\ref{qquasi}) into (\ref{niff_tra_g}) obtaining:
\begin{align}
&\left((G_{\sigma}(x))_i - (\overline G_\sigma(x))_i\right)^2 \leq  \nonumber\\
& L^2 \sigma^2 \, \int_{R^{n-1}} \int_{-\infty}^{\infty} (s_i^2 + \|\bar{s}_i\|^2)\, s_i^4\,\varphi(s_i)\, d s_i \> \> \varphi(\bar{s_i}) \, d {\bar{s_i}} + \nonumber\\
& + L^2 \sigma^2 \, \int_{R^{n-1}} \int_{-\infty}^{\infty} (s_i^2 + \|\bar{s}_i\|^2)\,
\bar{s_i}^2\,\varphi(s_i)\, d s_i \> \> \varphi(\bar{s_i}) \, d {\bar{s_i}}.\label{24}
\end{align}
For the first term in (\ref{24}) we obtain that
\begin{align}
&\int_{R^{n-1}} \int_{-\infty}^{\infty} (s_i^2 + \|\bar{s}_i\|^2)\, s_i^4\,\varphi(s_i)\, d s_i \> \> \varphi(\bar{s_i}) \, d {\bar{s_i}}  \nonumber\\
&= \int_{R^{n-1}} \int_{-\infty}^{\infty} s_i^6\,\varphi(s_i)\, d s_i \> \> \varphi(\bar{s_i}) \, d {\bar{s_i}} \nonumber\\
&\>\> + \int_{R^{n-1}} \int_{-\infty}^{\infty}  \|\bar{s}_i\|^2\, s_i^4\,\varphi(s_i)\, d s_i \> \> \varphi(\bar{s_i}) \, d {\bar{s_i}}\nonumber\\
&= 15 + 3(n-1). \label{primo}
\end{align}

By similar computations the second term in (\ref{24}) becomes
\begin{align}
&\int_{R^{n-1}} \int_{-\infty}^{\infty} (s_i^2 + \|\bar{s}_i\|^2)\,
\|\bar{s_i}\|^2\,\varphi(s_i)\, d s_i \> \> \varphi(\bar{s_i}) \, d {\bar{s_i}}\nonumber\\
&= \int_{R^{n-1}} \int_{-\infty}^{\infty} s_i^2\,
\|\bar{s_i}\|^2\,\varphi(s_i)\, d s_i \> \> \varphi(\bar{s_i}) \, d {\bar{s_i}}\nonumber\\
&\>\>+\int_{R^{n-1}} \int_{-\infty}^{\infty} \|\bar{s}_i\|^4
\,\varphi(s_i)\, d s_i \> \> \varphi(\bar{s_i}) \, d {\bar{s_i}}\nonumber\\
&= (n-1)+3(n-1) = 4(n-1).\label{secondo}
\end{align}
In (\ref{primo}) and (\ref{secondo}) we used the property \cite{cramer1999mathematical} that for a zero mean Gaussian $z$ with variance $\sigma^2$:
$$
E[z^{d}] =  \begin{cases}
({d}-1)!!\, \sigma^2 , & \text{for $d$ even } \\
0 , & \text{for ${d}$ odd }
\end{cases}
$$
where $(d-1)\,!! = (d-1)(d-3)\cdots 3\cdot 1$ \quad
and that for any $z \sim \mathcal{N}(0,I_{n-1})$
$$
\int_{R^{n-1}} \|z\|^2 \varphi(z) \, d z = \int_{R^{n-1}} \sum_{i=1}^{n-1} z_i^2 \varphi(z) \, d z = n-1.
$$
By substituting (\ref{primo}) and (\ref{secondo}) in (\ref{24}) we finally obtain that
\begin{equation*}
\left((G_{\sigma}(x))_i - (\overline G_\sigma(x))_i\right)^2 \leq L^2\,\sigma^2\,(15 + 3(n-1)+4(n-1)),
\end{equation*}
which, applied to all the entries, proves the theorem.
\bibliographystyle{spmpsci_unsrt}
\bibliography{bibliography}









\end{document}